\documentclass[12pt]{amsart}

\usepackage{tikz}
\usepackage{pgfplots}
\usepackage{filecontents}
\usepackage{amsmath, amssymb}
\usepackage[unicode]{hyperref}

\def\d{{\,\rm d}}
\newtheorem{Thm}{Theorem}
\newtheorem{claim}{Claim}

\newtheorem{Lem}{Lemma}

\begin{document}

\title{Nicolaas Govert de Bruijn, the enchanter of friable integers}

\author{Pieter Moree}
\address{Max-Planck-Institut f\"ur Mathematik,Vivatsgasse 7, D-53111 Bonn, Germany}
%\curraddr{}
\email{moree@mpim-bonn.mpg.de}
%\thanks{}

\date{}



\dedicatory{In memoriam: Nicolaas Govert (`Dick') de Bruijn (1918-2012) }
%{\def\thefootnote{}
%\footnote{{\it Mathematics Subject Classification (2000)}.
%11T22, 11B83}}

\begin{abstract}
\noindent N.G. de Bruijn carried out fundamental work on integers having only small prime factors 
and 
the Dickman-de Bruijn function that arises on computing the density of those integers. In this
he used his earlier work on linear functionals and differential-difference equations. We review the relevant work and also
some later improvements by others. 
\end{abstract}

\maketitle

\section{Introduction}
The number theoretical work of Nicolaas Govert (`Dick') de Bruijn (1918-2012) comes into two highly distinct flavours: combinatorial and analytical. In his 
combinatorial number theoretical work de Bruijn even did some cowork with one of the all-time greats in this
area: Paul Erd\H{o}s (6 joint papers!). Here we will only discuss de Bruijn's work in analytic number theory. This was done mainly 
in two periods: 1948-1953 and 1962-1966. In the second period de Bruijn revisited his earlier subjects. Some of this
work is joint with Jacobus Hendricus (`Jack') van Lint (1932-2004) \cite{As, vanL}.\\ 
\indent In {\it sieve theory} (see, e.g., \cite{FI, HaRi} or for 
a very brief introduction Section \ref{spezialschar} below) one is interested in estimating in terms of elementary functions the number of integers in
a prescribed sequence having a prescribed factorization. Usually one is interested
in the integers in the sequence that are primes and if one cannot handle those, integers that have only few distinct prime factors. De Bruijn's work in analytic number theory belongs mainly to sieve theory, but
he is not restricting the number of prime factors of the integers that remain after sieving. This is a milder 
form of sieving and here
usually quite sharp estimates can be obtained. As a basic example one can take the friable number counting function
$\Psi(x,y)$. Let $P(n)$ denote the largest prime divisor
of an integer $n\geqslant 2$. Put $P(1)=1$. 
A number $n$ is said to be {\it $y$-friable}\footnote{Some authors use $y$-smooth. Friable
is an adjective meaning easily crumbled or broken.} if $P(n)\leqslant y$. We let $S(x,y)$ denote the set of integers $1\leqslant n\leqslant x$
such that $P(n)\leqslant y$. The cardinality of $S(x,y)$ is denoted by $\Psi(x,y)$. \\
\indent In 1930, Dickman\footnote{Karl Dickman (1861-1947) wrote his paper when
he was 69 years old, after retiring from the Swedish insurance business world.} \cite{Dickman} proved that
\begin{equation}
\label{dikkertje}
\lim_{x\rightarrow \infty}{\Psi(x,x^{1/u})\over x}=\rho(u),
\end{equation}
where the {\it Dickman function} (today often also called {\it Dickman-de Bruijn function}) $\rho(u)$
is defined by 
\begin{equation}
\label{defie}
\rho(u)=
\begin{cases}
1 & \text{for $0\leqslant u\leqslant 1$};\\
{1\over u}\int_{0}^1 \rho(u-t)\d t & \text{for $u>1$.}
\end{cases}
\end{equation}

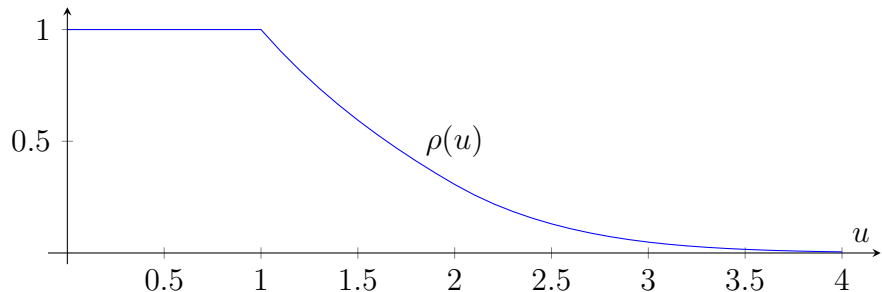
\begin{figure}
\caption{The Dickman-de Bruijn function $\rho(u)$}
\medskip
\begin{tikzpicture}
  \begin{axis}[
    width=\linewidth,
    height=5cm,
    axis x line=middle,
    axis y line=middle,
    xlabel=$u$, 
    xmin=-0.1, xmax=4.2,
    ymin=-0.05, ymax=1.1
    ]
    \addplot[color=blue] file {rho.dat};
    \node at (axis cs:2,0.5) {$\rho(u)$};
  \end{axis}
\end{tikzpicture}
\end{figure}
Note that $\rho(u)>0$, for if not, than because of the continuity of $\rho(u)$ there is a smallest zero $u_0>1$ 
and then substituting $u_0$ in (\ref{defie}) we easily arrive at a contradiction.
De Bruijn's contribution was to provide a very precise asymptotic estimate for $\rho(u)$ and in later works he provided
much more precise estimates then hitherto known for $\Psi(x,y)$. As a consequence various authors call $\rho(u)$ the
Dickman-de Bruijn function and some authors $\Psi(x,y)$ the Dickman-de Bruijn function. Both in depth and originality
de Bruijn's work goes way beyond everything done on the subject until then. Only starting in the 1980's, when Hildebrand and Tenenbaum 
started to work in this area, various results of de Bruijn were greatly improved upon. Most introductions of $\Psi(x,y)$ related
papers mention some de Bruijn papers. The techniques de Bruijn used in studying $\rho(u)$ are now standard in studying similar
functions that arise in sieve theory, see, e.g. \cite[Appendix B]{FI}.\\
\indent In the rest of this paper we first describe de Bruijn's work on $\rho(u)$ and then that on $\Psi(x,y)$. 
Before describing de Bruijn's work on $\rho(u)$ more in detail, we will discuss the basic properties of $\rho(u)$ and then
very briefly the linear functional work of de Bruijn.\\
\indent There is much more to friable integers than de Bruijn's work and its direct follow-up. The beautiful survey of Granville 
\cite{Granville}
makes very clear how many connections there are with other problems in (mainly) number theory, but also in probability, cycle structure of permutations (see Section \ref{avelar}), etc. etc. Here we will not go into this rich tapestry.\\
\indent It is assumed that the reader is familiar with the Landau-Bachmann O-notation (see wikipedia or any
introductory text on analytic number theory). Instead of $\log \log x$ we sometimes write $\log_2 x$, instead of
$(\log x)^A$, $\log^A x$.\\
\indent This paper is the extended version of a paper \cite{hollands} (in Dutch) written on request for
an issue of {\it Nieuw Archief voor Wiskunde} dedicated to the mathematics and memory of de Bruijn. For a discussion
of the work done on friable integers before de Bruijn, see Norton \cite{Norton} and Moree \cite{india}.

\section{Dick de Bruijn and the Dickman-de Bruijn function}
\subsection{Elementary properties of $\rho(u)$}
{}From (\ref{defie}) we see that $$\rho'(u)=-{\rho(u-1)\over u}$$ for 
$u>1$ and thus an alternative way of defining $\rho(u)$ is 
\begin{equation}
\label{defie2}
\rho(u)=
\begin{cases}
1 & \text{for $0\leqslant u\leqslant 1$};\\
1-\int_1^u {\rho(t-1)\over t}\d t & \text{for $u>1$.}
\end{cases}
\end{equation}
It follows that $\rho(u)=1-\log u$ for $1\leqslant u\leqslant 2$. 
For $2\leqslant u \leqslant 3$, $\rho(u)$ can be expressed in terms of the dilogarithm.
Put $\tau=(1+\sqrt{5})/2$. We have, see, e.g., Moree \cite{pietje},
$$\rho(\tau^2)=1-2\log \tau +(\log \tau)^2-{\pi^2\over 60}\approx 0.1046\ldots$$
The author does not know any other explicit pair $(u,v)$ satisfying $u>2$ and $\rho(u)=v$.\\
\indent Various papers 
discuss how to obtain {\it high accuracy numerical approximations} to $\rho(u)$, 
see, e.g., \cite{vdLW, mars}. The interest of de Bruijn was in the {\it asymptotical behaviour} of $\rho(u)$
and similar functions. So that will be our main focus from now on.\\
\indent Since $\rho'(u)=-\rho(u-1)/u$ for $u>1$ and $\rho(u)>0$ for
$u>0$, it follows that $\rho'(u)<0$ and hence $\rho(u)$ is monotonically decreasing for $u>1$. By (\ref{defie})
we find that $u\rho(u)\leqslant \rho(u-1)$, which on using
induction leads to
$\rho(u)\leqslant 1/[u]!$ for $u\geqslant 0$. It follows that $\rho(u)$ quickly tends to
zero as $u$ tends to infinity. De Bruijn was the first to obtain the following 
much more precise estimate for $\rho(u)$, for $u\geqslant 3$,
\begin{equation}
\label{eerstdick}
\rho(u)=\exp\Big\{-u\Big\{\log u+\log_2u-1+{\log_2 u-1\over \log u}+O\Big(\big({\log_2u\over \log u}\big)^2\Big)\Big\}\Big\}.
\end{equation}
Note that (\ref{eerstdick}) implies that, as $u\rightarrow \infty$,
$$\rho(u)={1\over u^{u+o(u)}},~~\rho(u)=\Big({{\rm e}+o(1)\over u\log u}\Big)^u,$$
formulae that suffice for most purposes and are easier to remember.\\
\indent De Bruijn obtained (\ref{eerstdick})  as a corollary to his 
highly non-elementary Theorem \ref{aso}, to be discussed below. In 
Section \ref{flauw}  an elementary proof of (\ref{eerstdick}) is given. 
\subsection{The Laplace transform of $\rho(u)$}
In de Bruijn's work the function $({\rm e}^x-1)/x$ and related functions play a prominent role. To understand its origin it is a very
useful exercise to compute the Laplace transform of $\rho(u)$. This transform is defined by
$${\hat \rho}(s)=\int_0^{\infty}\rho(u){\rm e}^{-us}\d u.$$
Using that $0<\rho(u)\leqslant 1/[u]!$ it is seen that this integral is absolutely convergent for all complex $s$, and
thus defines an entire fucntion of $s$.
Using (\ref{defie}) we then find that
\begin{align}
-{\hat \rho}'(s)&=\int_0^{\infty} u\rho(u){\rm e}^{-us}\d u =\int_0^1u{\rm e}^{-us}\d u+\int_1^{\infty}\Big(\int_{u-1}^u\rho(t)dt\Big){\rm e}^{-us}\d u\nonumber\\
&=\int_0^{\infty}\rho(t)\Big(\int_t^{t+1}{\rm e}^{-us}\d u\Big)\d t={1-{\rm e}^{-s}\over s}{\hat \rho}(s).\nonumber
\end{align}
It follows that ${\hat \rho}(s)=C\exp\{-{\rm Ein}(s)\}$, where $C$ is some constant and
$${\rm Ein}(s):=\int_0^s{1-{\rm e}^{-t}\over t}\d t,~s\in \Bbb C,$$
is the {\it complementary exponential integral}. The value of $C$ can be deduced by comparing the behaviour of 
${\hat \rho}(s)$ and $\exp\{-{\rm Ein}(s)\}$ as $s$ tends to infinity along the positive real axis. On the one hand we
have, by partial integration,
$$\int_0^{\infty}\rho'(u){\rm e}^{-us}\d u=-1+s{\hat \rho}(s)$$
and since the integral, in absolute value, is bounded by $1/s$, it follows that
$\lim_{s\rightarrow \infty}s{\hat \rho}(s)=1$. 
By integration by parts we see that
$$\int_0^1{{\rm e}^{-t}-1\over t}\d t+\int_1^{\infty}{{\rm e}^{-t}\over t}\d t=\int_0^{\infty}{\rm e}^{-t}\log t \d t=\Gamma'(1)=-\gamma,$$
where $\gamma$ denotes the Euler-Mascheroni constant and $\Gamma(z)$ the Gamma-function.
Thus, with $s$ tending to infinity along the positive real axis,
\begin{align}
{\rm Ein}(s)&=\int_0^1{1-{\rm e}^{-t}\over t}\d t-\int_1^s{{\rm e}^{-t}\over t}\d t+\log s\nonumber\\
&=\int_0^1{1-{\rm e}^{-t}\over t}\d t-\int_1^{\infty}{{\rm e}^{-t}\over t}\d t+\int_s^{\infty}{{\rm e}^{-t}\over t}\d t+\log s\nonumber\\
&=\gamma + \log s +O({\rm e}^{-s}).\nonumber
\end{align}
It follows that $C={\rm e}^{\gamma}$. Thus we proved (see also \cite[pp. 210-211]{monti} or \cite[pp. 370-372]{T})
 the following result. 
\begin{Lem}
\label{roe}
We have 
\begin{equation}
\label{fo}
{\hat \rho}(s)=\exp\big\{\gamma-\int_0^s{{1-{\rm e}^{-t}\over t}\d t}\big\}.
\end{equation}
\end{Lem}
As an amusing consequence it follows that for any $0\leqslant \delta\leqslant 1$,
$${\hat \rho}(0)={\rm e}^{\gamma}=\int_0^{\infty}\rho(t)\d t=\delta+\sum_{n\geqslant 1}(n+\delta)\rho(n+\delta).$$
On applying the inverse Laplace transform to both sides in (\ref{fo}) we obtain, for any real $\sigma_0$, 
\begin{equation}
\label{inverse}
\rho(u)={1\over 2\pi i}\int_{\sigma_0-i \infty}^{\sigma_0+i\infty}
\exp\Big\{\gamma+\int_0^s{{\rm e}^z-1\over z}dz+us\Big\}\d s,
\end{equation}
De Bruijn seems to have been the first to consider the Laplace transform of $\rho(u)$ (in \cite{B3}).
Lemma \ref{roe} is due to him.

\subsection{De Bruijn's $\xi$ function}
An important quantity in the asymptotic study of $\rho(u)$ is the function $\xi(u)$. For any given $u>1$, $\xi(u)$
is defined as the unique positive solution of the transcendental equation
\begin{equation}
\label{xixi}
{{\rm e}^{\xi}-1\over \xi}=u.
\end{equation}
It is easy to show that $\xi(u)$ is indeed unique. Put $f(x)=({\rm e}^x-1)/x$. Note that the Taylor series of $f(x)$ is
$$f(x)=1+{x\over 2!}+{x^2\over 3!}+{x^3\over 4!}+\ldots$$
So $f(x)$ tends to 1 as $x$ tends to $0+$ and is strictly increasing in $x$. This clearly shows that
$f(x)=u$ has a unique solution $x>0$ for each $u>1$.\\
\indent  {}From (\ref{xixi}) we infer, that as $u\rightarrow \infty$, we have
\begin{equation}
\label{xixixi}
\xi=\log \xi +\log u+O({1\over u\xi}).
\end{equation}
We want to find an expression for $\xi=\xi(u)$ for $u$ large. Note that for 
$u$ sufficiently large $1<\xi<2\log u$. Therefore $\log \xi =O(\log \log u)$. On feeding
this into the right hand side of (\ref{xixixi}) we obtain that $\xi=\log u +O(\log \log u)$, which on feeding again into 
the right hand side of (\ref{xixixi}) gives 
$$\xi=\log u + \log \log u+O\big({\log \log u\over \log u}\Big).$$
Iterating often enough we can get an error term
$O(\log^{-k} u)$ for any fixed $k>0$.\\
\indent By means of standard techniques of asymptotic analysis (see \cite[2.3]{Bboek}) it is possible to go beyond
this process of iteration. Using these a convergent series for $\xi(u)$ can be given (see \cite{HT2}).\\
\indent We will also need some information on $\int_1^u \xi(t)\d t$. By integration of a precise enough estimate
for $\xi(u)$ one obtains
\begin{equation}
\label{info}
\int_1^u \xi(t)\d t=u\Big\{\log u+\log_2u-1+{\log_2 u-1\over \log u}+O\Big(\big({\log_2u\over \log u}\big)^2\Big)\Big\}.
\end{equation}
As we will see the similarity between (\ref{info}) and (\ref{eerstdick}) is no coincidence (see Theorem \ref{A}).\\
\indent We leave it as an exercise to the reader to show that $$\lim_{u\rightarrow \infty}u\xi'(u)=1.$$

\subsection{Ramanujan's hand...}
On June 1-5, 1987, a meeting was held at the University of Illinois to commemorate the
centenary of Ramanujan's birth and to survey the many areas of mathematics (and of physics) 
that have been profoundly influenced by his work. One of the speakers, R. W. Gosper, remarked
in his lecture, "How can you love this man? He continually reaches his hand from his grave to snatch
your theorems from you." A few striking examples are given in an article by Berndt \cite{bruce}. 
One example involves Dickman's result (\ref{dikkertje}). Ramanujan \cite[p. 337]{Ramanujan} made the 
following claim 
(written at least ten years before Dickman's paper appeared!) that is
equivalent to the following claim.
\begin{claim}
Define, for $u\geqslant 0$, 
$$I_0(u)=1,~I_k(u)=
   \idotsint\limits_{t_1,\ldots,t_k\geqslant 1\atop t_1+\cdots+t_k\leqslant u}\cdots\int {\d t_1\cdots \d t_k\over t_1\cdots t_k},~k\geqslant 1$$
Then, for $u\geqslant 0$, {\rm (\ref{dikkertje})} holds with
$$\rho(u)=\sum_{k=0}^{\infty}{(-1)^k\over k!}I_k(u).$$
\end{claim}
For more details we refer to the book by Andrews and Berndt \cite[Section 8.2]{lost}. The asymptotic
behaviour of the $I_k(u)$ has been studied by Soundararajan \cite{sound}.
Related formulae for $\rho(u)$ as iterated integral can be, e.g., found in 
Buchstab \cite{Buchstab}, Chowla and Vijayaraghavan \cite{Chowla}, Goncharov \cite{go1} (see also
p. 22 of \cite{ABTbook}) and Tenenbaum \cite[(12) with $\lambda=1/u$ and $t=1$]{Told}.\\
\indent Ramanujan in his first letter (January 16th, 1913) to Hardy gave an asymptotic for $\Psi(x,3)$.
See Moree \cite{india} for further details.

\subsection{Mahler's partition problem}
\label{part}
Mahler's partition problem is to find an asymptotic formula for the
number $p_r(n)$ of representations of an integer $n$ in the form $n=n_0+n_1r+n_2r^2+\ldots$, $r>1$ a prescribed integer, in 
nonnegative integers $n_0,n_1,n_2,\ldots$, when $n$ is large. 
By considering the functional equation
\begin{equation}
\label{kurt}
F'(y)=F(y/r),
\end{equation}
K. Mahler \cite{maler, maler2} obtained an approximate formula for the number $p_r(n)$.
De Bruijn \cite{Dmal} determined the asymptotic behaviour much more precisely. A generalization of
de Bruijn's results was established by Pennington \cite{penning}, with $1,r,r^2,\ldots$ replaced by a more general
sequence $\lambda_1,\lambda_2,\ldots$ satisfying a certain growth condition.\\
\indent By setting $y=r^x$ and $F(r^x)=G(x)$, Mahler's functional equation is transformed into a functional equation
of the form
\begin{equation}
\label{kurtie}
G'(x)={\rm e}^{\alpha x+\beta}G(x-1),
\end{equation}
with $\beta$ real. In \cite{B4}, de Bruijn studied, using
saddle-point techniques, the asymptotic behavior of real solutions of the above equation as 
$x$ tends to infinity in case $\alpha>0$ and $\beta$ is a complex number.

\subsection{The linear functional equations papers}
\label{functional}
De Bruijn's initial mathematical interest was in combinatorics. According to J. Korevaar, de
Bruijn was really excited about
Mahler's partition problem (de Bruijn attended lectures of Mahler on this problem). This problem leads one to consider
linear functional equations and the author would not be surprised if this problem is the ultimate source of 
de Bruijn's strong
interest in these type of equations. Also he had a Ph.D. student who wrote
his Ph.D. thesis on this area of research \cite{been}.\\   
\indent In \cite{B-1}, de Bruijn studies the equation
$x^{-a}f′(x)+f(x)-f(x-1)=0$ and shows that depending on $a$ the behaviour of $f(x)$ is quite different. He
discusses an application to Mahler's partition problem (see Section \ref{part}).
In a later paper \cite{B0}, de Bruijn treats the 
more general equation $w(x)f′(x)+f(x)-f(x-1)=0$. It is shown that if $w(x)$ is not too small as $x$ tends
to infinity (a sufficient condition is $w(x)\geqslant 1/\log x$ as $x\rightarrow \infty$), then every solution tends to a constant $f(\infty)$, and an estimate is obtained for 
the difference $|f(x)-f(\infty)|$.\\ 
\indent In \cite{B2}, de Bruijn is concerned with equations of the type
$$f(x)=\int_0^1K(x,t)f(x-t)\d t.$$
He is interested in imposing conditions on $K$ which guarantee that every continuous solution $f(x)$ will be convergent, i.e., 
$\lim_{x\rightarrow \infty}f(x)$ exists and is finite (in such a case the kernel $K(x,t)$ is said to be `stabilizing').
As an application he establishes the following result. 
\begin{Thm}
\label{een}
If $F$ is continuous and satisfies 
\begin{equation}
\label{xF}
xF(x)=\int_0^1 F(x-t)\d t~(x>1),
\end{equation}
then there exists a 
constant $C$ such
that $F(x)=\{C+O(x^{-1/2})\}\rho(x)$.
\end{Thm}
De Bruijn used Theorem \ref{een} in order to prove Theorem \ref{aso}.

\subsection{De Bruijn's asymptotic formula for $\rho(u)$}
\label{dick1951}
In 1951, de Bruijn \cite{B3} proved the following theorem.
\begin{Thm} 
\label{aso}
As $u\rightarrow\infty$ we have
$$\rho(u)\sim {1\over \sqrt{2\pi u}}\exp\Big\{\gamma-u\xi(u)+\int_0^{\xi(u)}{{\rm e}^s-1\over s}\d s\Big\}.$$
\end{Thm}
Note that the integral on the right hand side equals $-{\rm Ein}(-\xi)$. We remark that
\begin{equation}
\label{E}
u\xi+{\rm Ein}(-\xi)=\int_1^u \xi(t)\d t,
\end{equation}
which can be seen on noting that
$${d\over du}(u\xi+{\rm Ein}(-\xi))=\big(u-{{\rm e}^{\xi}-1\over \xi}\big)\xi'(u)+\xi(u)=\xi(u),~u>1,$$
by the definition of $\xi$, and
$$\lim_{u\rightarrow 1+}(u\xi+{\rm Ein}(-\xi))=0,$$
since $\xi(u)\rightarrow 0$ as $u\rightarrow 1+$. Theorem \ref{aso} together with (\ref{E}) and
(\ref{info}) yield the estimate (\ref{eerstdick}) for $\rho(u)$.\\
\indent Below an outline of the remarkable proof of de Bruijn of Theorem \ref{aso}.\\
1) Note that
$$F_1(u)={1\over 2\pi i}\int_W \exp\Big\{-uz+\int_0^z{{\rm e}^s-1\over s}\d s\Big\}dz$$
satisfies (\ref{xF}), where the contour W consists of $\{z:\Re z\geqslant 0,~\Im z=-\pi\}$; followed by
$\{z:\Re z=0,-\pi\leqslant \Im z \leqslant +\pi\}$; followed by $\{z:\Re z\geqslant 0,~\Im z=+\pi\}$.\\
2) Use the saddle point method to estimate $F_1(u)$ for large $u$.\\
3) Note that for any continuous $G$ satisfying the adjoint equation
\begin{equation}
\label{dual}
uG'(u-1)=G(u)-G(u-1),~u>1,
\end{equation}
and any continuous $F$ satisfying $uF'(u-1)=-F(u-1),~u>1$, the function
$$\langle F,G\rangle =\int_{b-1}^b F(u)G(u)\d u - b F(b)G(b-1)$$
is a constant independent of $b>0$.\\
4) Note that
$$G_1(u)=\lim_{\delta\rightarrow 0}\Big\{\int_{-\infty}^{\delta}+\int_{\delta}^{\infty}\Big\}\exp\Big\{uz-\int_0^z{{\rm e}^s-1\over s}\d s\Big\}{{\rm e}^z\over z}dz,~u>-1,$$
satisfies (\ref{dual}).\\
5) Estimate $G_1(u)$, $u$ large, by the saddle point method.\\
6) Use the estimates of steps 2 and 5 to show that $\langle F_1,G_1\rangle=1$.\\
7) By considering appropriate limits, prove that $\langle \rho,G_1\rangle={
m e}^{\gamma}$.\\
8) By Theorem \ref{een} we then infer that $\rho(u)\sim {\rm e}^{\gamma}F_1(u)$.\\
9) Combine the latter result with the estimate obtained in step 2 in order to complete the proof. \qed\\

\noindent Applying the saddle-point method directly to (\ref{inverse}) it is possible to give a rather shorter
proof of Theorem \ref{aso}, see \cite[pp. 374--376]{T}\footnote{The first to remark this was Tenenbaum in 1990 in the French
original of his book.}. Indeed, in this way a rather sharper version of Theorem \ref{aso} can be obtained, see
Smida \cite{Smida}. It is not clear why de Bruijn did not choose to follow this easier approach.\\
\indent In 1982, Canfield \cite{Canfield} gave a reproof of Theorem \ref{aso} using combinatorial analysis. He defines two step functions
$p_N(u)$ and $q_N(u)$. The stepsize for each is $N^{-1}$; the value of $p_N$ on $[jN^{-1},(j+1)N^{-1})$
is $A(N,j)$ for $j\geqslant 0$; the value of $q_N$ on $((j-1)N^{-1},jN^{-1}]$ is $B(N,j)$ for $j\geqslant 1$, where
\begin{align}
{j\over N}A(N,j)&={1\over N}(A(N,j-1)+\ldots+A(N,j-N))\nonumber\\
A(N,0)&=A(N,1)=\ldots=A(N,N-1)=1;
\end{align}
and
\begin{align}
{j\over N}B(N,j)&={1\over N}(B(N,j-1)+\ldots+B(N,j-N+1))\nonumber\\
B(N,1)&=B(N,2)=\ldots=B(N,N-1)=1.
\end{align}
The idea is that the step functions $p_N$ and $q_N$ are discrete models of (\ref{defie}). The generating series
for $A(N,j)$ and $B(N,j)$ are not difficult to find:
$$\sum_{j=1}^{\infty}B(N+1,j)x^j=x\sum_{j=1}^{\infty}A(N,j)x^j=x\exp\big\{x+{x^2\over 2}+\ldots+{x^N\over N}\big\}.$$
Canfield then continues by showing that, for all $N$ and $u$,
$$B(N,[uN])\leqslant \rho(u)\leqslant A(N,[uN]).$$
Using the circle method he then shows there exists an integer-valued function $\nu(u)$ such that
$$B(\nu(u),[u\nu(u)])\sim A(\nu(u),[u\nu(u)])\sim 
{1\over \sqrt{2\pi u}}\exp(\gamma-u\xi-{\rm Ein}(-\xi)),$$
thus giving another proof of Theorem \ref{aso}.\\

\indent Alladi \cite{Alla} generalized de Bruijn's result as follows (see also \cite[p. 374]{T}).
\begin{Thm}
\label{A}
For $u\geqslant 1$ we have
$$\rho(u)=\Big(1+O({1\over u})\Big)\sqrt{\xi'(u)\over 2\pi}\exp\big(\gamma-\int_1^u\xi(t)\d t\Big).$$
\end{Thm}
Since $\lim_{u\rightarrow \infty}u\xi'(u)=1$ and 
because of the identity (\ref{E}), Theorem \ref{A} has Theorem \ref{aso} as a corollary.
Alladi writes that "this improvement is in fact obtained by iterating de Bruijn's method a second time, more 
carefully!". He starts by showing that in de Bruijn's Theorem \ref{een} the error term can be substantially improved.\\
\indent Hensley \cite{Hensley2} gave another proof of Theorem \ref{A}. An arithmetic proof of Theorem \ref{A} via the 
function $\Psi(x,y)$ is contained in Hildebrand and Tenenbaum \cite{HT0}. Hildebrand \cite{H2} established an asymptotic 
estimate similar to the one given in Theorem \ref{A} for every 
function in a one parameter family of differential-difference equations. Smida \cite{Smida} 
and Xuan \cite{Xuan} gave (quite technical) further
sharpenings of Theorem \ref{A}.\\
\indent De Bruijn in his proof of Theorem \ref{aso} apparently was the first to use an adjoint differential-difference
equation. Its application to questions arising in sieve theory was developed in several papers by
H. Iwaniec and others, in particular in Iwaniec \cite{rosser}.

\subsection{An elementary proof of de Bruijn's estimate {\rm (\ref{eerstdick})}}
\label{flauw}
Evertse et al. \cite[pp. 109-110]{EMST} showed that for $u\geqslant 1$ 
\begin{equation}
\label{me!}
\exp\Big(-\int_2^{u+1}\xi(t)\d t\Big)\leqslant \rho(u) \leqslant \exp\Big(-\int_1^{u}\xi(t)\d t\Big).
\end{equation}
The short and elementary proof of (\ref{me!}) only uses that $\rho'(u)/\rho(u)$ is non-decreasing for $u>1$. There is also an elementary
argument for that, see, e.g., \cite[pp. 35]{Moree}. Using the elementary estimate (\ref{info}) we infer from
(\ref{me!}) that
$$\rho(u)=\exp\Big(-\int_1^{u}\xi(t)\d t+O(\log u)\Big).$$
On inserting (\ref{info}) in this, the proof of (\ref{eerstdick}) is completed.\qed

\section{Preliminaries on $\Psi(x,y)$}
\subsection{Counting primes}
In order to study $\Psi(x,y)$ it is essential to have some understanding of the distribution of prime numbers. As usual
we let $\pi(x)=\sum_{p\leqslant x}1$ denote the number of primes $p\leqslant x$. We let li$(x)=\int_2^x \d t/\log t$ denote the
{\it logarithmic integral}. 
Hadamard and independently de la Vall\'ee Poussin in 1896 established the celebrated prime number
theorem, stating that 
\begin{equation}
\label{hadimassa}
\pi(x)\sim {x\over \log x},~x\rightarrow \infty.
\end{equation}
It follows from the proof of Hadamard and de la Vall\'ee Poussin that
\begin{equation}
\label{landau}
\pi(x)={\rm li}(x)+O(E(x)),
\end{equation}
with $$E(x)=x{\rm e}^{-c_1\sqrt{\log x}},$$
and $c_1$ some postive constant.
The currently best error term, due to I.M. Vinogradov and Korobov in 1958, is
\begin{equation}
\label{korobov}
E(x)=x\exp\{-c_2 \log^{3/5} x(\log_2 x)^{-1/5}\},
\end{equation}
with $c_2$ some positive constant (see \cite[Ch. 3]{Chandra}).
Important to keep in mind is that often
\begin{equation}
\label{approxi} 
\sum_{p\leqslant x}f(p)\approx \int_2^x {f(t)\over \log t}\d t.
\end{equation}
E.g., we expect that $\sum_{p\leqslant x}{1\over p}\approx \int_2^x {\d t\over t\log t}\approx \log \log x$.
Indeed, it is a classical result that
\begin{equation}
\label{klassiek}
\sum_{p\leqslant x}{1\over p}=\log \log x + c_3 +O({1\over \log x}),
\end{equation}
with $c_3$ a constant.\\
\indent Writing $\pi(x)={\rm li}(x)+E(x)$ we obtain, using the Stieltjes integral,
$$\sum_{p\leqslant x}f(p)=\int_2^x f(t)d\pi(t)=\int_2^x {f(t)\over \log t}\d t+\int_2^x f(t) \d E(t),$$
making (\ref{approxi}) more precise.
\subsection{The Buchstab functional equation for $\Psi(x,y)$}
\label{buchi}
Historically the first functional equation to be derived for the $\Psi(x,y)$ function was the
{\it Buchstab functional equation} \cite{Buchstab} (a few authors call this the Buchstab-de Bruijn equation),
\begin{equation}
\label{boek}
\Psi(x,y)=\Psi(x,z)-\sum_{y<p\leqslant z}\Psi({x\over p},p),
\end{equation}
where $1\leqslant y<z\leqslant x$. In the special case when $y=1$ and $z>1$ we obtain
$$\Psi(x,z)=1+\sum_{p\leqslant z}\Psi({x\over p},p).$$
In order to derive the Buchstab functional equation we start by observing that the difference $\Psi(x,z)-\Psi(x,y)$
equals precisely the cardinality of the set of natural numbers $\leqslant x$ having greatest prime factor in $(y,z]$. Let $p$
be some prime number. Notice that the cardinality of the set of natural numbers $\leqslant x$ having greatest prime factor
$p$ equals $\Psi(x/p,p)$. Thus, the Buchstab functional equation results. It can be used to determine $\Psi(x,y)$ successively
in the regions $y\geqslant x$, $\sqrt{x}<y\leqslant x$, $x^{1/3}\leqslant y<\sqrt{x}$. In terms of $u:=\log x/\log y$ this corresponds with
$u\leqslant 1,~1<u\leqslant 2,2<u\leqslant 3\ldots$. Note that $u$ is very close to the quotient of the number of digits of $x$ and the
number of digits of $y$.\\
\indent For $u\leqslant 1$ we obviously have $\Psi(x,y)=[x]$ and so in particular $\Psi(x/p,p)=[x/p]$ for $\sqrt{x}\leqslant p\leqslant x$. On taking $z=x$ in the Buchstab functional equation it then follows, for $1\leqslant u\leqslant 2$, that
\begin{equation}
\label{simpeltjes}
\Psi(x,y)=[x]-\sum_{y<p\leqslant x}[{x\over p}].
\end{equation}
Proceeding in this way one gets an exact expression for $\Psi(x,y)$ in increasingly large $u$ regions. Suppose we
have computed $\Psi(x,y)$ up to $u\leqslant h$. For a term $\Psi(x/p,p)$ in the Buchstab functional equation the logarithm
of the first argument  divided by the logarithm of the second argument equals $(\log x/\log p)-1$ and this is $<u-1$. So once
we can compute $\Psi(x,y)$ up to $u\leqslant h$ we can compute $\Psi(x,y)$ up to $u\leqslant h+1$. With each step of this process the
expression for $\Psi(x,y)$ becomes more complicated. The idea now is to invoke analytic number theory to `smoothen' the
sums. {}From (\ref{simpeltjes}) and (\ref{hadimassa}) we obtain
\begin{equation}
\label{simpeltjes2}
\Psi(x,y)=x-\sum_{y<p\leqslant x}{x\over p}+O(\pi(x))=x-\sum_{y<p\leqslant x}{x\over p}+O({x\over \log x}).
\end{equation}
This together with (\ref{klassiek}) yields, for all $x\geqslant 2$ and $1<u\leqslant 2$,
$$\Psi(x,x^{1/u})=x(1-\log u)\big(1+O({1\over \log x})\big).$$
This suggests that $\Psi(x,x^{1/u})\sim x f(u)$ for some function $f$ for every $u>0$. Heuristically this can be readily
justified. In the Buchstab functional equation we put $y=x^{1/u}$ and $z=x^{1/v}$ with $u>v>1$ arbitrary real numbers. We
then infer that
\begin{align}
xf(u) &\approx xf(v)-\sum_{x^{1/u}<n\leqslant x^{1/u}}{x\over p}f\big({\log x\over \log p}-1\big)\nonumber\\
& \approx x f(v)-\int_{x^{1/u}}^{x^{1/v}}{x\over t}f\big({\log x\over \log t}-1\big){\d t\over \log t}\nonumber\\
&= xf(v)- x \int_{v}^{u}{f(w-1)\over w}dw,\nonumber
\end{align}
where in the second step we used the approximation principle (\ref{approxi}) and in the final step the transformation
$t=x^{1/w}$ was made. Assuming that the involved error terms on dividing by $x$  tend to zero as $x$ and $y$ tend
to infinity, we conclude that the function $f(u)$ satisfies
$$f(u)=f(v)-\int_v^u{f(w-1)\over w}dw.$$
Thus $f(u)=\rho(u)$.  

\subsection{The Hildebrand functional equation for $\Psi(x,y)$}
The {\it Hildebrand identity} \cite{Hzeta} reads
$$\Psi(x,y)\log x=\int_1^x{\Psi(t,y)\over t}\d t+\sum_{p^m\leqslant x\atop p\leqslant y}\Psi({x\over p^m},y)\log p.$$
The advantage over the Buchstab identity is that one of the parameters is held fixed and 
that the coefficients are non-negative. It was the use of
this identity that led to serious improvements of several of de Bruijn's results. For most purposes it is enough
to work with the approximate identity
$$\Psi(x,y)\log x =\sum_{p\leqslant y}\Psi({x\over p},y) \log p+O(E(x,y)),$$
where one can take $E(x,y)=x$. For all $y\geqslant \log^{2+\epsilon}x$ one can take $E(x,y)=\Psi(x,y)$. An heuristic proof
of $\Psi(x,y)\sim x\rho(u)$ based on the Hildebrand functional equation can also be given, see \cite[pp. 416-417]{HT1}.
\subsection{Rankin's method}
In his work in the 30's on the gaps between consecutive primes Rankin \cite{Rankin} introduced a simple
idea to estimate $\Psi(x,y)$ which turns out to be remarkably effective and can be used in similar situations. Starting
point is the observation that for any $\sigma>0$
\begin{equation}
\label{rankinne}
\Psi(x,y)\leqslant \sum_{n\in S(x,y)}({x\over n})^{\sigma}\leqslant x^{\sigma}\sum_{P(n)\leqslant y}{1\over n^{\sigma}}=x^{\sigma}\zeta(\sigma,y),
\end{equation}
where 
$$\zeta(s,y)=\prod_{p\leqslant y}(1-p^{-s})^{-1},$$
is the partial Euler product up to $y$ for the Riemann zeta function $\zeta(s)$. Recall that, for $\Re s>1$,
$$\zeta(s)=\sum_{n=1}^{\infty}{1\over n^s}=\prod_{p}{1\over 1-p^{-s}}.$$
Note that the logarithmic derivative of $x^s\zeta(s,y)$ with respect to $s$ equals
$$\log x-\sum_{p\leqslant y}{\log p\over p^s-1}.$$
It is minimized for $\sigma=\alpha(x,y)$, where $\alpha(x,y)$ is the (unique) solution of
$$\sum_{p\leqslant y}{\log p\over p^{\alpha(x,y)}-1}=\log x.$$
{}From a sufficiently sharp form of the prime number theorem one, see \cite{HT0},
can derive the estimate
$$\alpha(x,y)={\log(1+y/\log x)\over \log y}\Big(1+O\Big({\log \log(1+y)\over \log y}\Big)\Big).$$
The upper bound $x^{\alpha}\zeta(\alpha,y)$ is only too large asymptotically by a 
factor $$\alpha \sqrt{2\pi\phi_2(\alpha,y)},$$ where $\phi_k(s,y)$ is defined as the $k$th partial
derivative of $\log \zeta(s,y)$ with respect to $s$. 
This factor turns out to be relatively small.
This follows from one of the deepest results in this area, due to Hildebrand and 
Tenenbaum \cite{HT0}, namely
that uniformly in the range $x\geqslant y\geqslant 2$ one has
$$\Psi(x,y)={x^{\alpha}\zeta(\alpha,y)\over \alpha \sqrt{2\pi\phi_2(\alpha,y)}}\Big(1+O({1\over u})+O\Big({\log y\over y}\Big)\Big).$$
\indent Of course, in (\ref{rankinne}) one is free to make any choice of $\sigma$. E.g., the choice
$\sigma=1-1/(2 \log y)$ leads to
$$\zeta(\sigma,y)\ll \exp \Big\{\sum_{p\leqslant y}{1\over p^{\sigma}}\Big\}\leqslant \exp\Big\{
\sum_{p\leqslant y}{1\over p}+O\Big((1-\sigma)\sum_{p\leqslant y}{\log p\over p}\Big)\Big\}\ll \log y,$$
which gives rise to $\Psi(x,y)\ll x{\rm e}^{-u/2}\log y$. The factor 
$\log y$ can be easily removed, see \cite[p. 359]{T}.\\ 
\indent As a further example of working with Rankin's method, let us try to estimate $\Psi(x,\log^A x)$ for $A>1$.
Letting $\sigma=1-1/A$, we get
$$\log \zeta(\sigma ,y)\ll \sum_{p\leqslant y}p^{-\sigma}=\sum_{p\leqslant y}{p^{1/A}\over p}\ll {y^{1/A}\over 
\log y}\ll {\log x\over \log \log x}.$$
We conclude that 
\begin{equation}
\label{ongelijk}
\Psi(x,\log^A x)\leqslant x^{1-1/A+O(1/\log \log x)}.
\end{equation}
We will see shortly that the upper bound is actually sharp.
\subsection{A binomial lower bound for $\Psi(x,y)$}
Let $2=p_1<p_2<p_3<\ldots$ be the consecutive primes. Let $p_k$ be the largest prime $p\leqslant y$ (thus $k=\pi(y)$).
Evidently $n$ is in $S(x,y)$ if and only if we can write $n$ in the form $n=p_1^{e_1}p_2^{e_2}\cdots p_k^{e_k}$, where
the $e_k$ are non-negative integers and $n\leqslant x$, that is,
\begin{equation}
\label{lineair}
e_1\log p_1+e_2\log p_2+ \cdots+e_k \log p_k\leqslant \log x.
\end{equation}
Since each $\log p_i\leqslant \log y$, any solution to
$e_1+e_2+\cdots+e_k\leqslant [{\log x\over \log y}]=[u]$ gives a solution to (\ref{lineair}) with $x=y^u$.
It follows that
\begin{equation}
\label{binom}
\Psi(x,x^{1/u})\geqslant \sum_{e_1+\ldots+e_k\leqslant [u]\atop e_i\geqslant 0}1= {[u]+\pi(y)\choose \pi(y)}.
\end{equation}
(If the reader cannot prove the last equality, (s)he can consult \cite[p. 208]{monti}.)

\section{De Bruijn's work on $\Psi(x,y)$}
\subsection{De Bruijn's $\Lambda$ function}
Let us define (as many authors do) $\rho(u)=0$ for $u<0$.
In \cite{most} de Bruijn introduced the function $\Lambda(x,y)$. He defines it as follows:
\begin{equation}
\label{lambda}
\Lambda(x,y)=x\int_0^{\infty}\rho\Big({\log x-\log t\over \log y}\Big)d({[t]\over t}).
\end{equation}
Since $\rho(u)=0$ for $u<0$ and $\rho(u)$ is continuous everywhere apart from the jump at $u=0$, the
integral is well-defined, unless the jump of $\rho$ coincides with a jump of $[t]/t$. This is the case if
$x$ is a positive integer $n$; this turns out to be a minor nuisance and can be dealt with by defining
$\Lambda(n,y)=\Lambda(n+0,y)$. The function $\Lambda(x,y)$ satisfies an equation which is similar to
the Buchstab equation (\ref{boek}). It arises from (\ref{boek}) by writing the sum in that equation formally
as a Stieltjes integral $\int \Psi(x/t,t)d\pi(t)$ and replacing $\pi(t)$ by li$(t)$. One can check that indeed
\begin{equation}
\label{boek2}
\Lambda(x,y)=\Lambda(x,z)-\int_{y}^{z}\Lambda({x\over t},t){\d t\over \log t}.
\end{equation}
Further, for $0<x\leqslant y$ we have $$\Lambda(x,y)=x\int_0^xd({[t]\over t})=x{[x]\over x}=[x]=\Psi(x,y).$$ Thus
$\Lambda(x,y)$ satisfies the same initial condition as $\Psi(x,y)$  and obeys a functional equation which is
the friable version of Buchstab's identity.\\
\indent De Bruijn showed that $\Lambda(x,y)$ closely approximates $\Psi(x,y)$ in a certain
$(x,y)$-region (this region was substantially extended in 1989 by Saias \cite{Saias}).
Further he shows that both
$\Lambda(x,y)$ and $\Psi(x,y)$ satisfy an expansion of the form
$$x\sum_{r=0}^n a_r{\rho^{(r)}(u)\over \log^r y}+E_n(x),$$
with $E_n(x)$ a term he describes rather explicitly and
$$K(z)=a_0+a_1z+a_2z^2+\ldots$$ the Taylor series around $z=0$ of $z\zeta(1+z)/(1+z)$. 
The term $E_n(x)$ is only an error term if $u$ is not close to an integer from above. 
For a deeper understanding of this phenomenon, see
Theorem 1.1 and Lemma 3.1 in Hanrot et al.~\cite{HTW}.

\subsection{De Bruijn's uniform version of Dickman's result}
In \cite{most} de Bruijn gives estimates for the differences 
$$|\Psi(x,y)-\Lambda(x,y)|{\rm ~and~}|\Lambda(x,y)-x\rho(u)|.$$
These estimates involve the error term $E(x)$ in the prime number theorem, cf. (\ref{landau}).
These bounds, when combined with the sharpest known form of the prime
number theorem (having error term (\ref{korobov})), yield the following uniform estimate of Dickman's 
limit result (\ref{dikkertje}).
\begin{Thm}
\label{vijf}
The estimate 
\begin{equation}
\label{bruni}
\Psi(x,y)=x\rho(u)\Big\{1+O\Big({\log(u+1)\over \log y}\Big)\Big\},
\end{equation}
holds for $1\leqslant u\leqslant \log^{3/5-\epsilon}y,{~that~is,~}y>\exp(\log^{5/8+\epsilon}x)$.
\end{Thm}
Hildebrand \cite{H} improved this substantially to the range
$$1\leqslant u\leqslant \exp(\log^{3/5-\epsilon}y),{\rm ~that~is,~}y>\exp((\log \log x)^{5/3+\epsilon}).$$
There are limitations as to how far this range can be made. Hildebrand \cite{Hzeta} showed
that the estimate (\ref{bruni}) holds uniformly for
$$1\leqslant u\leqslant y^{1/2-\epsilon},{\rm ~that~is,~}y\geqslant \log^{2+\epsilon}x,$$
if and only if the Riemann Hypothesis is true.\\
\indent One might wonder in which region $\Psi(x,y)\gg x \rho(u)$ is valid, as did Hensley \cite{Hensley}.
Canfield et al.~\cite{CEP} proved that 
there is a constant $C$ such that
$\Psi(x,y)\geqslant x\rho(u)\exp\{C(\log_2u/\log u)^2\}$ uniformly for $x\ge 1$ and $u\ge 3$.
It seems that $\Psi(x,y)\geqslant x\rho(u)$ for a very large $(x,y)$-region.

\subsection{The average largest prime factor of $n$}
\label{avelar}
Recall that $P(n)$ denotes the largest prime factor of $n$.
De Bruijn \cite{most} applies his results on $\Lambda(x,y)$ to show
that
\begin{equation}
\label{logpee}
\sum_{n\leqslant x}\log P(n)=\lambda x \log x+O(x),{\rm ~with~}
\lambda=\int_0^{\infty}{\rho(u)\over (1+u)^2}\d u.
\end{equation}
By partial integration it follows that
$$\lim_{x\rightarrow \infty}{1\over x}\sum_{2\leqslant n\leqslant x}{\log P(n)\over \log n}=\lambda,$$
thus proving an heuristic claim by Dickman \cite{Dickman}. The interpretation of the result is that
for an average integer with $m$ digits, its greatest prime factor has about $\lambda m$ digits.
De Bruijn's estimate can be further sharpened to
$$\sum_{n\leqslant x}\log P(n)=\lambda x \log x-\lambda(1-\gamma)x+O(x\exp(-\log^{3/8-\epsilon}x)),$$
see \cite[Exercise III.5.3]{T}. 
The constant $\lambda$ is now known as
the {\it Golomb-Dickman constant}. 
Mitchell \cite{Mitchell} computed that $\lambda=0.62432 99885 4\dots$.
Golomb et al.~\cite{golomb} have defined a constant $\mu$ which is related to
the limiting properties of random permutations. Let $L_n$ be the expected length of the longest cycle of a random
permutation of $n$ letters. Define $\mu_n=L_n/n$. (Thus $\mu_1=1$, $\mu_2=3/4$, $\mu_3=13/18$, $\mu_4=67/96$.)
It can be shown that the numbers $\mu_n$ are monotonically decreasing with $n$ and thus a limit $\mu$ exists. It turns
out that $\lambda=\mu$. Thus on average the longest cycle of a permutation of $n$ elements has 
length $\lambda n$. 
Much earlier Goncharov announced in \cite{go1} and proved in \cite{go2} that the ratio of permutations in $S_n$ having
longest cycle $\leqslant n/u$, with $u\ge 1$ and fixed, tends to $\rho(u)$ as $n$ tends to infinity. Goncharov, unaware
of the occurrence of $\rho(u)$ in number theory,
also expressed $\rho(u)$ as an iterated integral, see 
Arratia et al. \cite[p. 22]{ABTbook}.
Shepp and Lloyd  \cite{SL} derived a second formula for $\lambda$:
$$\lambda=\int_0^{\infty}\exp\{-x-\int_x^{\infty}{{\rm e}^{-y}\over y}\d y\}\d x.$$
It is now known that the cycle lengths in a permutation are typically Poisson \cite{grandeur}.\\
\indent It is striking that $\rho(u)$ turns up both in the setting of greatest cycle length and greatest prime factor. It also
shows up in the setting of greatest irreducible polynomial factor of a polynomial over a finite field, see, e.g., 
Car \cite{car}. For a unified treatment see Arratia et al. \cite{ABT} (or their book \cite{ABTbook}).\\
\indent Meanwhile there are many results (see, e.g., \cite{EIP,scour}) involving $P(n)$, one of the more interesting, see \cite{EIP}, being
$$\sum_{n\leqslant x}{1\over P(n)}=x\int_2^x\rho({\log x\over \log t}){\d t\over t^2}\Big(1+O\Big(\big({\log_2 x\over \log x}\big)^{1/2}\Big)\Big).$$
Later Scourfield \cite{scour} improved on the error term.
Sums of the form
$\sum_{n\leqslant x}f(P(n))$, where $f$ is some arithmetic function can be evaluated using the identity
$$\sum_{n\leqslant x}f(P(n))=\sum_{p\leqslant x}f(p)\Psi({x\over p},p),$$
together with sufficiently sharp estimates for $\Psi(x,y)$.\\
\indent A quite general result is due to Tenenbaum and Wu \cite{TW4},
who gave an asymptotic expansion for
$\sum_{n\leqslant x} f(n)\{\log P(n)\}^r$, where $f$ is a fairly general non-negative arithmetical function (which need not be multiplicative) and $r$ is any positive real number.

\subsection{De Bruijn's 1966 estimate for $\log \Psi(x,y)$}
We now turn our attention to de Bruijn's 1966 paper \cite{B5}. In it he is concerned with
estimating $\log \Psi(x,y)$.
Set 
\begin{equation}
\label{Z}
Z:={\log x\over \log y}\log\Big(1+{y\over \log x}\Big)+{y\over \log y}\log\Big(1+{\log x\over y}\Big).
\end{equation}
\begin{Thm} 
\label{1966}
We have, uniformly for $x\geqslant y\geqslant 2$,
$$\log \Psi(x,y)=Z\Big\{1+O\Big({1\over \log y}+{1\over u}+{1\over \log \log 2x}\Big)\Big\}.$$
\end{Thm}
Actually, in the error term $u^{-1}$ can be left out, see Tenenbaum \cite[pp. 359--362]{T} for a proof.
Let us apply Theorem \ref{1966} with $y=\log^A x$ and $A>1$ fixed. One then finds that 
$$\Psi(x,\log^A x)=x^{1-1/A+O(1/\log \log x)},$$
thus improving on (\ref{ongelijk}). Already the results of de Bruijn in \cite{most} imply as a special case that
$\Psi(x,\log^A x)=x^{1-1/A+O(\epsilon)}$. Following de Bruijn's argument very closely Chowla and Briggs \cite{ChowlaB} gave
a reproof that they believed to be considerably simpler.\\
\indent The proof of Theorem \ref{1966} is very elementary (but takes several pages
to work out). The upper bound is deduced from Rankin's upper bound (\ref{rankinne})
with an optimal choice of the parameter $\sigma$. The lower bound is based from the inequality (\ref{binom})  together
with a lower bound for the right hand side of (\ref{binom})  obtained by Stirling's formula.\\
\indent Theorem \ref{1966} clearly shows
that there is a change of behaviour of $\Psi(x,y)$ as $y\approx \log x$. If $y/\log x\rightarrow \infty$, then the
first term in (\ref{Z}) dominates, whereas $Z$ is asymptotic to the second term in
(\ref{Z}) when $y=o(\log x)$. The change in behaviour is due to the fact that if $y$ is small compared to $\log x$, then many prime
factors of a `typical' integer in $S(x,y)$ occur to high powers, whereas for larger values of $y$ most prime factors occur only
to the first power. 
P. Erd\H{o}s \cite{paultje} showed that, if $x\rightarrow \infty$, we have
$$\log \Psi(x,\log x)\sim {(\log 4)\log x\over \log \log x}.$$

\section{De Bruijn and the sieve of Eratosthenes}
Define $\Phi(x,y)$ as the number of integers $n\leqslant x$ having no prime factors $\leqslant y$.
In the sieve of Eratosthenes all the integers $1\leqslant n\leqslant x$, which are multiples of the primes $2\leqslant p<\sqrt{x}$ are
removed. What remains after this process is the number 1 and all the prime numbers $p$ in the range $\sqrt{x}\leqslant p\leqslant x$.
It follows that  
\begin{equation}
\label{sieve}
\Phi(x,\sqrt{x})=1+\pi(x)-\pi(\sqrt{x}).
\end{equation}  
Buchstab \cite{Buchstab} showed that
\begin{equation}
\label{boek1}
\Phi(x,y)\sim x{\rm e}^{\gamma}\omega(u)\prod_{p<y}(1-{1\over p}),
\end{equation}
where $\omega(u)=1/u$ for $1\leqslant u\leqslant 2$ and, for all $u\geqslant 2$,
$$u\omega(u)=1+\int_1^{u-1}\omega(t)\d t.$$
By a theorem of Mertens we have $\prod_{p<y}(1-1/p)^{-1}\sim {
m e}^{-\gamma}/\log y$, and thus we can reformulate (\ref{boek1})
as
$$\Phi(x,y)\sim u\omega(u){x\over \log x}.$$
It follows, e.g., that $\Phi(x,\sqrt{x})\sim x/\log x$, an asymptotic one also finds from (\ref{hadimassa})
and (\ref{sieve}).\\
\indent De Bruijn \cite{B1} writes
$$\Phi(x,y)=x\prod_{p<y}(1-p^{-1})\Psi_1(x,y)$$
and shows that
$$|\Psi_1(x,y)-{\rm e}^{\gamma} \log y \int_1^u y^{t-u}\omega(t)\d t|<c_1\exp(-c_2\sqrt{\log y}),~u\geqslant 1,~y\geqslant 2,$$
for appropriate constants $c_1$ and $c_2$. \\
\indent De Bruijn \cite[Example 1]{B0} showed that $\lim_{u\rightarrow \infty}\omega(u)$ exists, and, if we denote it by $A$, then
\begin{equation}
\label{zegeensA}
\omega(u)=A+O({1\over [u]!}). 
\end{equation}
In a later paper de Bruijn \cite{B1} showed that $A={\rm e}^{-\gamma}$. De Bruijn stated that the right order
of magnitude of $|\omega(u)-{\rm e}^{-\gamma}|$ "is probably something like 
\begin{equation}
\label{diic}
\exp (-u\log u-u \log_2 u){\rm ~}".
\end{equation}  
This was proved
subsequently by L.K. Hua \cite{Hua} by an "advanced calculus argument" and by 
Buchstab \cite{Buchstab2} using an arithmetic argument. The formula (\ref{diic}) reminds one of the 
asymptotic estimate (\ref{eerstdick}) and indeed it turns out to be natural to write $\omega(u)-e^{-\gamma}=\rho(u)h(u)$ in
order to give a sharp estimate for the difference, see Tenenbaum \cite[Lemme 4]{kriebel}.
A further connection between $\omega(u)$ and $\rho(u)$ is that their Laplace transforms 
${\hat \rho}(s)$ and ${\hat \omega}(s)$ are related by
$$1+{\hat \omega}(s)={1\over s{\hat \rho}(s)},~s\ne 0.$$
\indent A surprising and important application of the properties of $\omega(u)$ has been mady by Maier \cite{Maier}. By an
ingenious construction, Maier proved that the number of primes in short intervals $[x,x+\log^C x]$, $C>1$, is
sometimes larger than, and sometimes smaller than the expected number $\log^{C-1}x$. Put $W(u)={\rm e}^{\gamma}\omega(u)-1$.
Put $M_{+}(v)=\max_{u\geqslant v}W(u)$ and $M_{-}(v)=\min_{u\geqslant v}W(u)$. Maier proved for any fixed $C>1$,
$$\lim \sup_{x\rightarrow \infty}{\pi(x+\log^C x)-\pi(x)\over \log^{C-1}x}\geqslant 1+M_{+}(C);$$
$$\lim \inf_{x\rightarrow \infty}{\pi(x+\log^C x)-\pi(x)\over \log^{C-1}x}\leqslant 1+M_{-}(C);$$
Using a method due to de Bruijn involving the adjoint equation of $\omega(u)$, Maier showed that $W(u)$ changes
sign in every interval of length one. Hence, for all $C\geqslant 1$, $M_{+}(C)>0$ and $M_{-}(C)<0$. The paper
of Maier led to a lot of follow-up work on the lack of equi-distribution of the primes. 
In this work methods from the study of friable integers and differential-difference equations play an
important role, see e.g. \cite{FG, FGHM}.

\section{How special are $\rho$ and $\omega$?}
\label{spezialschar}
The reader might wonder how special $\rho(u)$ and $\omega(u)$ are. The answer is that many similar functions
occur in number theory, especially in sieve theory (the sieve of Eratosthenes being the easiest example of a sieve).
Let us give some examples.\\

\noindent 1) Let $Q$ be a set of primes, and denote by $\Psi(x,y,Q)$ the number of positive integers 
not exceeding $x$ that have no prime factors from $Q$ exceeding $y$. Suppose that the number of primes 
$\pi(x,Q)$ of primes $p\leqslant x$ in $Q$ satisfies
$$\pi(x,Q)=\delta {\rm li}(x)+O({x\over \log^{B} x}),$$
with $0<\delta<1$.
Then Goldston and McCurley \cite{gold} showed that
$$\Psi(x,y,Q)=x\tau_{\delta}(u)\Big(1+O\big({1\over \log y}\big)\Big),$$
uniformly for $u\geqslant 1$ and $y\geqslant 1.5$. 
Here $\tau_{\delta}(u)$ is a function similar to $\rho(u)$, it is the 
unique solution of $\tau_{\delta}(u)=1$ for $0\leqslant u\leqslant 1$
and $$u\tau'_{\delta}(u)=-\delta \tau_{\delta}(u-1){\rm ~for~}u>1.$$ This family of differential-difference equations
had been earlier investigated Beenakker \cite{been}, a Ph.D. student of de Bruijn.\\
\indent The latter result is not valid if $Q$ is the set of all primes (then 
$\delta=1$), as the $\log(u+1)$ in
the error term in (\ref{bruni}) is sharp and cannot be removed.\\

\noindent 2) Let $k\geqslant 1$ be an integer. Let $P_k(n)$ denote the $k$th largest prime factor of $n$. 
Let $\Psi_k(x,y)$ denote the number of integers $n\leqslant x$ such that $P_k(n)\leqslant y$.
Knuth and Trabb Pardo \cite{KP} showed that
$$
\lim_{x\rightarrow \infty}{\Psi_k(x,x^{1/u})\over x}=\rho_{1,k}(u),
$$
where $\rho_{1,k}(u)$ is a function similar to the Dickman-de Bruijn function.\\

\noindent 3) In sieving we consider a finite sequence $A$ of integers and a set $P$ of primes. Let $X$ denote the number of 
integers in $A$. The goal is to count the number of elements $S(A,P,z)$ that remain after sifting $A$ by all the
primes in $P$ less than $z$, that is we are interested in estimating
$$S(A,P,z):=|\{a\in A:(a,P(z))=1\}|,$$
where $P(z):=\prod_{p\in P,~p<z}p$. This formulation is much too general if one wants to obtain non-trivial results
and so we impose regularity conditions for $d|P(z)$ on the subsequence
$$A_d:=\{a\in A:a\equiv 0({\rm mod~}d)\}.$$
We assume that we have approximations of the form
\begin{equation}
\label{zeef1}
|A_d|={\omega(d)\over d}X+r_d,
\end{equation}
where $\omega(p)$ is a multiplicative function satisfying $0\le \omega(p)<p$ such that $r_d$ is a remainder that
is small. For $p\not\in P$ we put $\omega(p)=0$. We define
$$V(z):=\prod_{p<z}\big(1-{\omega(p)\over p}\big).$$
Intuitively $V(z)$ can be regarded as the probability that an element of $A$ is not divisible by any prime $p<z$ with
$p$ in $P$.  We further assume there is a positive constant $c_1$ such that 
\begin{equation}
\label{zeef2}
0\le {\omega(p)\over p}\le 1-{1\over c_1}
\end{equation}
for every $p$ in $P$ and that there are constants $\kappa>0$ and $c_2\ge 1$ such that
\begin{equation}
\label{zeef3}
\Big| \sum_{w_1\le p<w_2}{\omega(p)\over p}\log p-\kappa \log({w_2\over w_1})\Big|\le c_2,{\rm ~for~all~}w_2\ge w_1\ge 2.
\end{equation}
The constant $\kappa$ is called the {\it dimension of the sieve} and is roughly the average of $\omega(p)$ over the primes 
as 
$$\sum_{p\le w_2}{\log p\over p}\sim \log w_2$$
as $w_2$ tends to infinity. Under the assumptions (\ref{zeef1}), (\ref{zeef2}) and (\ref{zeef3}) one obtains the following
inequalities for $S(A,P,z)$:
\begin{equation}
\label{zeef4}
S(A,P,z)\le XV(z)\{F_{\kappa}\big({\log y\over \log z}\big)+o\}+R_y,
\end{equation} 
and
\begin{equation}
\label{zeef5}
S(A,P,z)\ge XV(z)\{f_{\kappa}\big({\log y\over \log z}\big)+o\}+R_y,
\end{equation}
respectively, where $F_{\kappa}(u)$ and $f_{\kappa}(u)$ are {\it sieve auxiliary functions} depending on the
dimension $\kappa$, $y$ is a free parameter, $o$ is an error term and $R_y$ is an error term formed of the $r_d$'s of
(\ref{zeef1}), cf. the inequality (\ref{zeef6}). In applications the free parameter $y$ is chosen in such a way that the
remainder term $R_y$ is small in comparison with the main term containing $X$.\\
\indent The Buchstab functional equation (\ref{boek}) has a counterpart in sieve theory:
$$S(A,P,z)=S(A,P,z_1)-\sum_{z_1\le p<z,~p\in P}S(A_p,P,p),~z\ge z_1\ge 2.$$
{}From this we can obtain a lower bound for $S(A,P,z)$ on assuming that a lower bound for $S(A,P,z_1)$ is already known, and
that we have an upper bound for the terms in the sum. For $z_1$ small enough an estimate for $S(A,P,z)$ follows from the
so called Fundamental Lemma of Sieve theory which states that under the regularity conditions (\ref{zeef1}), (\ref{zeef2})
and (\ref{zeef3}) supplemented with the condition that $|r_d|\le \omega(d)$ for $d|P(z)$ and $d$ square free, we have
for $z_1\le X$,
$$S(A,P,z_1)=XV(z_1)\{1+O({\rm e}^{-u/2})\},$$
where $u=\log X/\log z_1$ \cite[Theorem 2.5]{HaRi}. \\
\indent As a demonstration of the above let us consider the Selberg's upper bound sieve 
\cite{HaRi, raw}. Let $\kappa\ge 0$, $u=\log y/\log z>0$ and $\nu(d)$ the number of distinct prime divisors of $d$.
Under the conditions (\ref{zeef1}), (\ref{zeef2}), (\ref{zeef3}), we have
\begin{equation}
\label{zeef6}
S(A,P,z)\le {XV(z)\over \sigma_{\kappa}(u)}\big\{1+O\Big({(\log_2 y)^{2\kappa+1}\over \log y}\Big)\Big\}+
\sum_{d<y,~d|P(z)}3^{\nu(d)}|r_d|,
\end{equation}
where $\sigma_{\kappa}(u)$ is the Ankeny-Onishi-Selberg function. The implied constant depends only on $\kappa,c_1$ and $c_2$.
One has $\sigma_{\kappa}(u)=\rho_{\kappa}(u)u^{\kappa-1}$, where $\rho_{\kappa}(u)=1$ for $0\le u\le 1$ and
$$\rho_{\kappa}'(u)=-u^{-\kappa}(1-u)^{\kappa-1}\rho_{\kappa}(u-1),~u>1.$$
Note that $\sigma_1(u)=\rho(u)$.
It seems that $\sigma_{\kappa}(u)$ was introduced independently by Ankeny and Onishi \cite{AO} and
de Bruijn and van Lint \cite{BL}.\\
\indent The function $\sigma_{\kappa}(u)$ satisfies the differential-difference equation
\begin{equation}
\label{verg}
uf'(u)+af(u)+bf(u-1)=0,
\end{equation}
with $a=1-\kappa$ and $b=\kappa$. The Buchstab sifting function $\omega(u)$ satisfies (\ref{verg}) with $a=1$ and
$b=-1$. Solutions to the equation (\ref{verg}) appear to have been first studied by Iwaniec \cite{rosser} in
connection with his work on Rosser's sieve. A more systematic investigation was carried out by
Wheeler \cite{wheeler}. Both Iwaniec and Wheeler make intensive use of the so called adjoint equation 
\begin{equation}
\label{adjointje1}
ug'(u)+(1-a)g(u)-bg(u+1)=0,
\end{equation}
which, in some sense, is easier to deal with than (\ref{verg}). Whereas the functions satisfying (\ref{verg}) can
exhibit a rather erratic behaviour, it turns out that there exists a solution of (\ref{adjointje1}) which is analytic
in the right half plane. Since solutions to the equations (\ref{verg}) and (\ref{adjointje1}) are connected to a simple
integral relation (the right hand side of (\ref{inner}), which turns out to be a constant), one can derive asymptotic
information for a given solution to (\ref{verg}) by studying the asymptotic behaviour of a suitable solution to
(\ref{adjointje1}). Notice that this is precisely the approach followed by de Bruijn in proving Theorem \ref{aso}! It
has been established by Hildebrand and Tenenbaum \cite{HT2} that also in this general case a much more direct
approach is possible (cf. the paragraph in Section \ref{dick1951} following the proof of de Bruijn).\\
\indent The problem of describing the general solution to (\ref{verg}) amounts to describing solutions $f(u)=f(u;\varphi)$ 
satisfying (\ref{verg}) for $u>u_0$ and $f(u)=\phi(u)$ for $u_0-1\leqslant u\leqslant u_0$, where $u_0$ is any positive
real number and $\phi(u)$ is any given continuous function on $[u_0-1,u_0]$. Given two functions $f$ and $g$ defined
on $[u_0-1,u_0]$ and $[u_0,u_0+1]$, respectively, we set
\begin{equation}
\label{inner}
\langle f,g\rangle = u_0f(u_0)g(u_0)-b\int_{u_0-1}^{u_0}f(u)g(u+1)\d u.
\end{equation}
If $f$ and $g$ are solutions of (\ref{verg}), respectively (\ref{adjointje1}), 
then this "scalar product" is independent of $u_0$.\\
\indent Hildebrand and Tenenbaum \cite{HT2} showed that there is a solution $F(u;a,b)$ and that there are solutions $F_n(u;a,b)$ for every integer
$n$ that form a basis of the solution space, that is any solution can be expressed in the form
$$f(u)=\alpha F(u;a,b)+\sum_{n=-\infty}^{\infty}\alpha_n F_n(u;a,b),$$
with suitable coefficients $\alpha$ and $\alpha_n$ depending on the initial function $\phi$. 
De Bruijn \cite{B4} by a different approach had earlier found a similar result for solutions of the differential-difference equation (\ref{kurtie}).\\
\indent Earlier special cases of (\ref{verg}) had been considered by Alladi \cite{al3}, Beenakker \cite{been}, 
Hensley \cite{Hensley2}, Hildebrand \cite{H2} and Wheeler \cite{wheeler} amongst others.

\section{Arithmetic sums over friable integers}
\label{poech}
In \cite{BL} de Bruijn and van Lint consider sums of the form $$\sum_{n\in S(x,y)}f(n),$$ with $f$ a 
non-negative multiplicative function. Many of the usual techniques can be applied here, e.g.,
the analogues of the Buchstab and Hildebrand functional equation hold, see e.g. \cite{morf}.\\
\indent For an enlightening historical introduction and an extensive bibliography on the subject, see
Tenenbaum and Wu \cite{TW1}. 
For the state of the
art of this area see Tenenbaum and Wu \cite{TW4}.\\
\indent Let $\Psi_m(x,y)$ denote the number of $y$-friables that are $\leqslant x$ and coprime to $m$.
It turns out that evaluating the ratio $\Psi_m(x/d,y)/\Psi(x,y)$ for $1\leqslant d\leqslant x,~m\geqslant 1,~x\geqslant y\geqslant 2$, is a crucial step for estimating arithmetic sums over friable integers, see de la Bret\`eche and
Tenenbaum \cite{delaB1}.
Alladi \cite{al1, al2} noticed a duality which shows that $\Psi(x,y)$ is related to the sum 
of the M\"obius function over the uncancelled elements in the sieve 
of Eratosthenes, and the sum of the M\"obius function over the smooth 
numbers is related to $\Phi(x,y)$.\\  
\indent An interesting related theme is that of the {\it friable Tur\'an-Kubilius inequality}. 
The Tur\'an-Kubilius inequality $$\sum_{n\leqslant x}|f(n)-A_f(x)|^2\ll xB_f(x),$$ holds uniformly for strongly additive functions $f$, where $A_f(x)$, $B_f(x)$ stand for certain sums over primes $p\leqslant x$, depending on $f$. 
If $$B_f(x)=o(A_f(x),$$ it can be shown that
$f$ has {\it normal order} $A_f(x)$, that is, for
any $\epsilon>0$, we have $$|f(n)-A_f(n)|\leqslant \epsilon|A_f(n)|$$ on a set of integers $n$ of density 1. 
One can try
to obtain an inequality analogous to the Tur\'an-Kubilius inequality for the corresponding 
sum over the elements of the set $S(x,y)$. The first result in this direction is due to
Alladi \cite{Alla}, who found that he needed a stronger version (Theorem \ref{A}) of 
de Bruijn's Theorem \ref{aso} for this purpose. For further developments the reader is referred to
works of Tenenbaum and his collaborators \cite{delaB1, delaB2, delaB3, HMT, mate}.

\section{De Bruijn's analytic number theory work in a nutshell}
\centerline{--\cite{Dmal}, 1948--}
\vskip 4mm
\noindent Mahler's partition problem is to find an asymptotic formula for $p_r(n)$, the number of 
partitions of $n$ into powers of a fixed integer $r>1$. De Bruijn gives an asymptotic for $\log p_r(n)$ improving
on Mahler's result \cite{maler}  (cf. Section \ref{part}).
\vskip 2mm
\centerline{--\cite{B-1}, 1949--}
\vskip 2mm
\noindent This paper is mostly in pure analysis. At the end of the paper there is a brief discussion of connections with
Mahler's partition problem and a promise to come back to this in a future paper (\cite{B4}).
\vskip 2mm
\centerline{--\cite{B0}, 1950--}
\vskip 2mm
\noindent As an application of more general results de Bruijn shows that $\omega(u)=A+O(1/[u]!)$ for some constant $A$.
\vskip 2mm
\centerline{--\cite{B1}, 1950--}
\vskip 2mm
\noindent He gives a sharp estimate for $\Phi(x,y)$.
Further he gives the first proof that the Buchstab sifting function $\omega(u)$ tends to ${\rm e}^{-\gamma}$ as $u$ tends to infinity 
and formulates a conjecture (see (\ref{diic}) that was later proved 
by `the father of modern Chinese mathematics' (Loo-Keng Hua) \cite{Hua}.
\vskip 2mm
\vfil\eject
\centerline{--\cite{B2}, 1950--}
\vskip 2mm
\noindent He proves Theorem \ref{een}.
\vskip 2mm
\centerline{--\cite{B3}, 1951--}
\vskip 2mm
\noindent The main result he proves here is Theorem \ref{aso}, further
Lemma \ref{roe} and the inverse Laplace transform formula (\ref{inverse}). 
As a corollary of his main theorem he obtains (\ref{eerstdick}).
Further, he introduces the function $\xi(u)$.
\vskip 2mm
\centerline{--\cite{most}, 1951--}
\vskip 2mm
\noindent He introduces $\Lambda(x,y)$, shows that it closely approximates both $\Psi(x,y)$ and $x\rho(u)$. These estimates in combination with the 
error estimate (\ref{korobov}) in the prime number theorem lead to a uniform version
of Dickman's result (\ref{dikkertje}), namely Theorem \ref{vijf}. Further, he uses his $\Lambda(x,y)$
function to prove the estimate (\ref{logpee}).
\vskip 2mm
\centerline{--\cite{B4}, 1953--}
\vskip 2mm
\noindent He studies, using saddle-point techniques, the asymptotic behavior 
in $x$ of real solutions of the differential-difference
equation $$F'(x)={\rm e}^{\alpha x+\beta}F(x-1),~{\rm with~}\alpha>0,\beta\in \mathbb C.$$ 
This is related to 
Mahler's partition problem (cf. Section \ref{part}), in which case $\beta$ is real.
\vskip 2mm
\centerline{--\cite{B4a}, 1962/1963--}
\vskip 2mm
\noindent If $n$ and $x$ are positive integers, let $f(n,x)$ denote the number of integers $1\le m\le x$ such that
all prime factors of $m$ divide $n$. Let $\gamma(k)=\prod_{p|k}p$ denote the squarefree kernel of $k$. Erd\H{o}s in a letter to
de Bruijn conjectured that $F(x):=\sum_{n\leqslant x}f(n,x)$ satisfies $\log(F(x)/x)=O(\log^{1/2+\epsilon}x)$. Put $G(x)=\sum_{n\leqslant x}f(n,n)$. De Bruijn, using a Tauberian theorem by Hardy and
Ramanujan, proves that
$$\log({F(x)\over x})\sim \log({G(x)\over x})\sim \log \big(\sum_{n\leqslant x}{1\over \gamma(n)}\big)\sim \sqrt{8\log x\over \log_2 x},$$
thus proving a stronger form of Erd\H{o}s conjecture. (Later W. Schwarz \cite{zwart} gave an asymptotic for
$\sum_{n\leqslant x}\gamma(n)^{-1}$ itself.)
\vskip 2mm
\centerline{--\cite{B5}, 1966--}
\vskip 2mm 
\noindent De Bruijn establishes Theorem \ref{1966}, using Rankin's method, a binomial lower bound
for $\Psi(x,y)$ and Stirling's formula.
\vskip 2mm
\vfil\eject
\centerline{--\cite{BL-1}, 1963--}
\vskip 2mm 
\noindent De Bruijn and van Lint show that $F(x)\sim G(x)$ (in the notation of \cite{B4a}), thus answering a question
raised by Erd\H{o}s.
\vskip 2mm
\centerline{--\cite{BL}, 1964--}
\vskip 2mm 
\noindent De Bruijn and van Lint consider sums of the form $\sum_{n\in S(x,y)}f(n)$, with $f$ a 
non-negative multiplicative function (cf. Section \ref{poech}). They introduced a generalization
of the Dickman-de Bruijn function (see Section \ref{spezialschar}) that around the same time showed up in the analysis of Ankeny
and Onishi \cite{AO} of the Selberg upper bound sieve.
\subsection{De Bruijn's writing style}
When de Bruijn introduces a new result or method, he also describes to what extent this improves
on earlier work. 
Interesting special cases of his result he also likes to consider.
However, he also discusses limitations and might also describe
an alternative method. Thus he looks at a topic from various sides and gives an honest evaluation of
merits and drawbacks of his work. 
In this context, I was struck by the fact that in a lecture he gave
at age 90 \cite{90} he tried to give an honest evaluation of the positive and negative sides of his own character!
In terms of giving details of proofs he is a bit on the brief side (as
was usual in mathematical research papers written in 1945-1955). 
Often he promises to come back to some aspect in a future paper (and keeps his promise!).\\
\indent De Bruijn's book \cite{Bboek} because of its clarity and pleasant style found many adepts. Rather
than being theoretical the book proceeds through discussing interesting example, several of 
which are from combinatorics,
e.g. asymptotics for the Bell numbers\footnote{De Bruijn's book predates the time when
the Bell numbers had this name.}. De Bruijn's book
was the starting point and inspiration for several later works
specializing on the use of analytic methods in combinatorics and 
algorithms, e.g., \cite{Bender, Od, FlaS}\footnote{This paragraph follows closely an e-mail E.R. Canfield sent me.}.

\section{Further reading}
\label{further}
In the first place the reader is advised to consult the book \cite{Bboek} and the papers of de Bruijn!
As a first introduction to friable numbers I can highly recommend Granville's 2008 survey \cite{Granville}. 
It has a special
emphasis on friable numbers and their role in algorithms in computational number theory. 
Mathematically more demanding is
the 1993 survey by Hildebrand and Tenenbaum \cite{HT1}. 
The literature up to 1971 is discussed in extenso in
Norton's monograph \cite{Norton}. For a discussion of the work prior to 1950 (when de Bruijn entered the scene), see
Moree \cite{india}.
Chapter III.5 in Tenenbaum's book \cite{T}
deals with $\rho(u)$ and approximations to $\Psi(x,y)$ by the saddle point method, Chapter III.6 with the dual problem
of counting integers $n\leqslant x$ having no prime factors $\leqslant y$. For an introductory account of the saddle point method, 
see, e.g., Tenenbaum \cite{Tcol}.\\

\noindent {\tt Acknowledgement}. 
The surveys by Granville \cite{Granville} and Hildebrand and Tenenbaum \cite{HT1} were very helpful to me.
In Sections \ref{buchi} and \ref{spezialschar} I copied close to verbatim material (not published in a journal) from my 
Ph.D. thesis \cite{Moree}. 
I copied close to verbatim the outline of de Bruijn's proof of Theorem \ref{aso} in 
Section \ref{dick1951} from
Canfield \cite{Canfield}. 
The 
title of this paper is inspired by the title of a biography of Ada Lovelace \cite{ada}. 
I thank J. Sorenson and A. Weisse for their kind assistance in creating Figure 1.
A. Vershik pointed out to me the work of W. Goncharov \cite{go1, go2} and E. Bach kindly sent me \cite{go1}. 
K. Alladi, A.D. Barbour, D. Bradley, E.R. Canfield, L. Holst, A. Ivic, J. Korevaar, K.K. Norton, J. Sorenson, P. Tegelaar and R. Tijdeman are heartily thanked
for their comments. My special thanks are due to G. Tenenbaum for his very extensive comments and helpful
correspondence.

\end{document}